\documentclass{amsart}
\usepackage{amssymb}

\newtheorem{Lem}{Lemma}[section]
\newtheorem{Facts}{Facts}[section]
\newtheorem{lemma}[Lem]{Lemma}
\newtheorem{Thm}[Lem]{Theorem}
\newtheorem{Cor}[Lem]{Corollary}
\theoremstyle{definition}
\newtheorem{Def}[Lem]{Definition}
\newtheorem{Asm}[Lem]{Assumption}

\DeclareMathOperator{\dom}{dom}
\DeclareMathOperator{\ult}{Ult}
\DeclareMathOperator{\cf}{cf}
\newcommand{\ON}{\text{ON}}
\newcommand{\BM}{\text{BM}}
\newcommand{\PD}{\text{PD}}

\newcommand{\name}[1]{\dot{#1}}

\begin{document}

\date{2007-01-23}
\title{Winning the pressing down game but not Banach Mazur}

\author[Jakob Kellner]{Jakob Kellner$^\ast$}
\address{Kurt G\"odel Research Center for Mathematical Logic\\
  Universit\"at Wien\\
  W\"ahringer Stra\ss e 25\\
  1090 Wien, Austria}
\email{kellner@fsmat.at}
\urladdr{http://www.logic.univie.ac.at/$\sim$kellner}
\thanks{$^\ast$ supported by a European Union Marie Curie EIF fellowship, contract MEIF-CT-2006-024483.}

\author{Matti Pauna}
\address{Department of Mathematics and Statistics\\
  University of Helsinki\\
  Gustaf H\"allstr\"omin katu 2b\\
  FIN-00014,  Finland}
\email{matti.pauna@helsinki.fi}
\urladdr{http://www.helsinki.fi/$\sim$pauna/}

\author[Saharon Shelah]{Saharon Shelah$^\dag$}

\address{Einstein Institute of Mathematics\\
  Edmond J. Safra Campus, Givat Ram\\
  The Hebrew University of Jerusalem\\
  Jerusalem, 91904, Israel\\
  and
  Department of Mathematics\\
  Rutgers University\\
  New Brunswick, NJ 08854, USA}
\email{shelah@math.huji.ac.il}
\urladdr{http://shelah.logic.at/}
\thanks{$^\dag$ supported by
  the United States-Israel Binational Science Foundation (Grant no. 2002323),
  publication 896.}

\begin{abstract}
  Let $S$ be the set of those $\alpha\in\omega_2$ that have
  cofinality $\omega_1$.
  It is consistent relative to a measurable that
  the nonempty player wins
  the pressing down game of length $\omega_1$,
  but not 
  the Banach Mazur game of length $\omega+1$
  (both games starting with $S$).
\end{abstract}

\maketitle

\section{Introduction}

We set 
$E^{\kappa}_\theta=\{\alpha\in\kappa:\, \cf(\alpha)=\theta\}$.
Let $S$ be a stationary set.  We investigate two games, each played by players
called ``empty'' and ``nonempty''. Empty has the first move.

In the Banach Mazur game $\BM(S)$ of length $<\theta$, the players choose
decreasing stationary subsets of $S$. Empty wins, if at some $\alpha<\theta$
the intersection of these sets is nonstationary. (Exact definitions are give in
the next section.)

In the pressing down game $\PD(S)$, empty cannot choose a stationary subset of
the moves so far, but only a regressive function. Nonempty chooses a homogeneous
stationary subset.

So it is at least as hard for nonempty to win $\BM$ as to win $\PD$. 

In this paper, we show that $\BM$ can be really harder than $\PD$. 
This follows easily from well known facts about precipitous ideals
(cf. \ref{cor:withcfomega} for a more detailed explanation):
Nonempty can never win $\BM_{\leq\omega}(\omega_2)$, but it is consistent
(relative to a measurable) that nonempty wins $\PD_{<\omega_1}(\omega_2)$. The
reason is the following: In $\BM$, empty can first choose
$E^{\omega_2}_{\omega}$, and empty always wins on this set. However in $\PD$,
it is enough for nonempty to win on $E^{\omega_2}_{\omega_1}$.
In a certain way this is ``cheating'', since nonempty wins $\PD$ on
$E^{\omega_2}_{\omega_1}$ but looses $\BM$ on the disjoint set
$E^{\omega_2}_{\omega}$. So in a way the difference arises
because empty has the first move in $\BM$.

Therefore, a better question is:
Can nonempty win $\PD(S)$ but loose $\BM(S)$ even if nonempty gets the first
move,\footnote{Which is equivalent to: nonempty does not win
  $\BM_{\leq\omega}(S')$ for any stationary $S'\subseteq S$.}
e.g. on $S=E^{\omega_2}_{\omega_1}$?

This is indeed the case:
\begin{Thm}\label{thm:main}
   It is consistent relative to a measurable that
   for $\theta=\aleph_1$ and $S=E^{\theta^+}_\theta$,
   nonempty wins $\PD_{<\omega_1}(S)$ but not $\BM_{\leq\omega}(S)$,
   even if nonempty gets the first move.

   The same holds for
   $\theta=\aleph_n$ (for $n\in \omega$) etc.
\end{Thm}

Various aspects of these and related games have been studied for a long time.

Note that in this paper we consider the games on sets, i.e. a move is an
element of the powerset of $\kappa$ minus the (nonstationary) ideal. A popular
(closely related but not always equivalent) variant are games on a Boolean
algebra $B$: Moves are elements of $B$, in our case $B$ would be the 
powerset of $\kappa$ {\em modulo} the ideal. 

Also note that in Banach Mazur games of length greater than $\omega$, it is
relevant which player moves first at limit stages (in our definition this is
the empty player).  Of course it is also important who moves first at stage $0$
(in this paper again the empty player), but the difference here comes down to a
simple density effect (cf. \ref{lem:trivial}.4).

The Banach Mazur $\BM$ game has been investigated e.g. in
\cite{MR519808} or \cite{MR846604}.
It is closely related to the so-called ``ideal game'' and to precipitous
ideals, cf. Theorem \ref{thm:old} and \cite{MR0505504}, \cite{MR0485391}, or
\cite{MR560220}.
$\BM$ is also related to the ``cut \& choose game'' of \cite{MR739910}.

The pressing down game is related to the Ehrenfeucht-Fra\"iss\'e game in model
theory, cf. \cite{MR1191613} or \cite{MR1971240}, and has applications in set
theory as well \cite{MR1469093}.

Other related games have been studied e.g. in  \cite{MR0739911} or
\cite{MR1420262}.

We thank Jouko V{\"a}{\"a}n{\"a}nen for asking about Theorem \ref{thm:main} and
for pointing out Theorem \ref{cor:inner}.

\section{Banach Mazur, pressing down, and precipitous ideals}

Let $\kappa$ and $\theta$ be regular, $\theta<\kappa$.

We set
$E^\kappa_\theta=\{\alpha\in\kappa:\, \cf(\alpha)=\theta\}$.
$\mathcal{E}^\kappa_\theta$ is the family of stationary subsets
of $E^\kappa_\theta$.
Analogously
for $E^\kappa_{>\theta}$ etc.

Instead of ``the empty player has a winning strategy for the game $G$'' we just
say ``empty wins $G$'' (as opposed to: empty wins a specific run of the game).

$\mathcal{I}$ denotes a fine, normal ideal on $\kappa$ (which implies
$<\kappa$-completeness).

A set  $S\subseteq \kappa$ is called $\mathcal{I}$-positive if
$S\notin \mathcal{I}$.

\begin{Def} Let $\kappa$ be regular, and
  $S\subseteq \kappa$ an $\mathcal{I}$-positive set.
  \begin{itemize}
    \item
      $\BM_{< \zeta}(\mathcal{I},S)$, the
      {\em Banach Mazur game} of length
      $\zeta$ starting with $S$, is played as follows:

      At stage $0$, empty plays an $\mathcal{I}$-positive
      $S_0\subseteq S$, nonempty plays $T_0\subseteq S_0$.
      At stage $\alpha<\zeta$, empty plays an $\mathcal{I}$-positive
      $S_\alpha \subseteq \bigcap_{\beta<\alpha}S_\beta$ (if possible),
      and nonempty plays some  $T_\alpha\subseteq S_\alpha$.

      Empty wins the run, if 
      $\bigcap_{\beta<\alpha}S_\beta\in \mathcal{I}$
      at any stage $\alpha<\zeta$.  Otherwise nonempty wins.

      (For nonempty to win a run, it is not necessary that
      $\bigcap_{\beta<\zeta}S_\beta$ is $\mathcal{I}$-positive or even just nonempty.)

    \item $\BM_{\leq \omega}(\mathcal{I},S)$ is $\BM_{<\omega+1}(\mathcal{I},S)$. So 
      empty wins the run iff $\bigcap_{n<\omega}S_n\in \mathcal{I}$.

    \item $\PD_{< \zeta}(\mathcal{I},S)$, the {\em pressing down game}
      of length $\zeta$ starting with $S$, is played as follows:

      At stage $\alpha<\zeta$, empty plays a regressive function 
      $f_\alpha:\kappa\to\kappa$, and nonempty plays some $f_\alpha$-homogeneous
      $T_\alpha\subseteq \bigcap_{\beta<\alpha}T_\beta$.

      Empty wins the run, if $T_\alpha\in \mathcal{I}$
      for any $\alpha<\zeta$.
      Otherwise, nonempty wins.

    \item $\PD_{\leq \omega}(\mathcal{I},S)$ is $\PD_{<\omega+1}(\mathcal{I},S)$.

    \item $\BM_{<\zeta}(S)$ is $\BM_{<\zeta}(\text{NS},S)$, 
      and $\PD_{<\zeta}(S)$ is $\PD_{<\zeta}(\text{NS},S)$ (where NS denotes
      the nonstationary ideal).
  \end{itemize}
\end{Def}

The following is trivial:
\begin{Facts}\label{lem:trivial}
  \begin{enumerate}
    \item Assume $S\subseteq T$.
      \begin{itemize}
        \item
          If empty wins $\BM_{<\zeta}(\mathcal{I},S)$, then
          empty wins $\BM_{<\zeta}(\mathcal{I},T)$.
        \item
          If nonempty wins $\BM_{<\zeta}(\mathcal{I},T)$, then
          nonempty wins $\BM_{<\zeta}(\mathcal{I},S)$.
        \item
          If empty wins $\PD_{<\zeta}(\mathcal{I},T)$, then
          empty wins $\PD_{<\zeta}(\mathcal{I},S)$.
        \item
          If nonempty wins $\PD_{<\zeta}(\mathcal{I},S)$, then
          nonempty wins $\PD_{<\zeta}(\mathcal{I},T)$.
      \end{itemize}
    \item Assume $\mathcal{I}\subseteq \mathcal{J}$, and
      $\mathcal{J}$ also is fine and normal.
      \begin{itemize}
        \item
          If empty wins $\PD_{<\zeta}(\mathcal{I},S)$, then
          empty wins $\PD_{<\zeta}(\mathcal{J},S)$.
        \item
          If nonempty wins $\PD_{<\zeta}(\mathcal{J},S)$, then
          nonempty wins $\PD_{<\zeta}(\mathcal{I},S)$.
      \end{itemize}
    \item In particular, if nonempty wins $\PD_{<\zeta}(\mathcal{I},S)$, then nonempty wins
       $\PD_{<\zeta}(S)$.
    \item Let $\BM'$ be the variant of $\BM$ where nonempty gets the
      first move (at stage 0 only). The difference between $\BM$ and $\BM'$
      is a simple density effect:
      \begin{itemize}
        \item Empty wins $\BM'_{<\zeta}(\mathcal{I},S)$ iff
	  empty wins $\BM_{<\zeta}(\mathcal{I},S')$ for all positive
          $S'\subseteq S$ iff empty has a winning strategy for $\BM$ 
          with $S$ as first move.
        \item Empty wins $\BM_{<\zeta}(\mathcal{I},S)$ iff 
          empty wins $\BM'_{<\zeta}(\mathcal{I},S')$ for some positive 
          $S'\subseteq S$.
        \item Nonempty wins $\BM'_{<\zeta}(\mathcal{I},S)$ iff  nonempty wins
          $\BM_{<\zeta}(\mathcal{I},S')$ for some positive $S'\subseteq S$.
      \end{itemize}
  \end{enumerate}
\end{Facts}
(For 3, use that $\mathcal{I}$ is normal, which implies 
$\text{NS}\subseteq \mathcal{I}$.)

We will use the following definitions and facts concerning precipitous ideals,
as introduced by Jech and Prikry
\cite{MR0505504}.
We will usually refer to Jech's {\em Millennium Edition} \cite{MR1940513} for
details.

\begin{Def} Let $\mathcal{I}$ be a normal ideal on $\kappa$.
  \begin{itemize}
    \item Let $V$ be an inner model of $W$. 
      $U\in W$ is called a normal $V$-ultrafilter 
      if the following holds:
      \begin{itemize}
        \item If $A\in U$, then $A\in V$ and $A$ is a subset of $\kappa$. 
        \item $\emptyset\notin U$, and $\kappa\in U$.
        \item If $A,B\in V$ are subsets of $\kappa$, $A\subseteq B$ and $A\in U$, then $B\in U$.
	\item If $A\in V$ is a subset of $\kappa$, then either $A\in U$ or
          $\kappa\setminus A\in U$.
        \item If $f\in V$ is a regressive function on $A\in U$, then
          $f$ is constant on some $B\in U$.
      \end{itemize}
      (Note that we do not require iterability or amenability.)
    \item A normal $V$-ultrafilter $U$ is wellfounded, if the ultrapower of 
      $V$ modulo $U$ is wellfounded. In this case the transitive collapse
      of the ultrapower is denoted by $\ult_U(V)$.
    \item Let $P_\mathcal{I}$ be the  
      family of $\mathcal{I}$-positive sets
      ordered by inclusion.
      Since $\mathcal{I}$ is normal, 
      $P_\mathcal{I}$ forces that the generic filter $G$ is
      a normal $V$-ultrafilter (cf. \cite[22.13]{MR1940513}).
      $\mathcal{I}$ is called precipitous, if $P_\mathcal{I}$ forces that 
      $G$ is wellfounded.
    \item The ideal game on $\mathcal{I}$ is played just like
      $\BM_{\leq\omega}(\mathcal{I},\kappa)$,
      but empty wins iff $\bigcap_{n\in\omega} S_n$ is empty 
      (as opposed to ``in $\mathcal{I}$'').
  \end{itemize}
\end{Def}

So if empty wins the ideal game, then empty wins
$\BM_{\leq\omega}(\mathcal{I},\kappa)$. And if nonempty wins
$\BM_{\leq\omega}(\mathcal{I},\kappa)$, then nonempty wins the ideal game.

\begin{Thm}\label{thm:old} Let $\mathcal{I}$ be a normal
  ideal on $\kappa$.
  \begin{enumerate}
    \item (Jech, cf \cite[22.21]{MR1940513})
      $\mathcal{I}$ is not precipitous iff empty wins the ideal game. So in this
      case empty also wins $\BM_{\leq \omega}(\mathcal{I},\kappa)$.
    \item (cf. \cite{MR0485391})
      If $E^\kappa_\omega \notin \mathcal{I}$, then
      nonempty cannot win the ideal game, and empty wins\footnote{
        There is even a fixed sequence of winning moves for empty:
	For every $\alpha\in E^\kappa_\omega$ let
        $(\alpha_n)_{n\in\omega}$ be a normal sequence in $\alpha$. 
        As move $n$, empty plays the function that maps $\alpha$ to $\alpha_n$.
        If $\beta$ and $\beta'$ are both in $\bigcap_{n\in \omega}T_n$, then
        $\beta_n=\beta'_n$ for all $n$ and therefore $\beta=\beta'$.}
      $\PD_{\leq \omega}(\mathcal{I},E^\kappa_\omega)$
      and therefore also $\BM_{\leq \omega}(\mathcal{I},\kappa)$.
    \item (Jech, Prikry \cite{MR560220}, cf \cite[22.33]{MR1940513}) 
      If $\mathcal{I}$ is precipitous, then $\kappa$ is measurable in an inner
      model.
    \item (Laver, see \cite{MR0485391} or \cite[22.33]{MR1940513})
      Assume that $U$ is a normal ultrafilter on $\kappa$.
      Let $\aleph_1\leq \theta<\kappa$ be regular and let
      $Q = \text{Levy}(\theta, <\kappa)$ be the Levy collapse (cf lemma
      \ref{lem:basiclevy}).
      In $V[G_Q]$, let $\mathcal F$ be the filter generated by $U$ and $\mathcal{I}$ the
      corresponding ideal. Then $\mathcal{I}$ is normal, and the
      family of $\mathcal{I}$-positive sets has a $<\theta$-closed dense subfamily.

      So in particular it is forced that 
      nonempty wins $\BM_{<\theta}(\mathcal{I},S)$ for all
      $\mathcal{I}$-positive
      sets $S$ (nonempty just has to pick sets from the dense subfamily),
      and therefore that nonempty wins $\PD_{<\theta}(S)$ (cf \ref{lem:trivial}.3).
    \item (Magidor \cite{MR560220}, penultimate paragraph)
      One can modify this forcing to get
      a $<\theta$-closed dense subset 
      of $\mathcal{E}^{\theta^+}_\theta$.

      So in particular, $\mathcal E^{\theta^+}_\theta$ can be precipitous.
  \end{enumerate}
\end{Thm}

Mitchell \cite{MR560220} showed that the $\text{Levy}(\omega, <\kappa)$ gives a
precipitous ideal on $\omega_1$ (and with Magidor's extension, $\text{NS}_{\omega_1}$
can be made precipitous). So the ideal game is interesting on $\omega_1$,
but our games are not:

\begin{Cor}\label{cor:withcfomega}
  \begin{enumerate}
    \item  
      Empty always wins $\PD_{\leq\omega}(S)$
      (and $\BM_{\leq\omega}(S)$) for $S\subseteq \omega_1$.
    \item
      It is equiconsistent with a measurable
      that nonempty wins $\BM_{<\theta}(E^{\theta^+}_\theta)$
      for e.g. $\theta=\aleph_1$, $\theta=\aleph_2$,
      $\theta=\aleph^+_{\aleph_7}$ etc.
    \item
      The following is consistent relative to a measurable: Nonempty wins
      $\PD_{<\theta}(\theta^+)$ but not $\BM_{\leq \omega}(\theta^+)$
      for e.g. $\theta=\omega_1$.
  \end{enumerate}
\end{Cor}

\begin{proof}
  (1) is just \ref{thm:old}.2, and (2) follows from \ref{thm:old}.3--4.

  (3) Let $\kappa$ be measurable, and 
    Levy-collapse $\kappa$ to $\theta^+$.
    According to \ref{thm:old}.2, nonempty wins $\PD_{<\omega_1}(S)$
    for all $S\in U$, in particular for $S=\theta^+$.
    However, empty wins $\BM_{\leq \omega}(\theta^+)$ (by playing
    $E^{\theta^+}_\omega$).
\end{proof}

In the rest of the paper will deal with the proof of Theorem 
\ref{thm:main}.

\section{Overview of the proof}

We assume that $\kappa$ is measurable, and $\omega<\theta<\kappa$ regular.

\subsection*{Step 1}
We construct models $M$ satisfying:\\
\centerline{($\ast$) $\kappa$ is measurable and player empty wins $\BM_{\leq
       \omega}(S)$ for every stationary $S$.}

We present two constructions, showing that ($\ast$) is true in $L[U]$
as well as compatible with larger cardinals:
{\renewcommand{\labelenumi}{(\roman{enumi})}
  \begin{enumerate}
    \item The inner model $L[U]$, Section \ref{sec:inner}:\\
       Let $D$ be a normal ultrafilter on
       $\kappa$, and set $U=D\cap L[D]$.
       Then in $L[U]$, (the dual ideal of) $U$ is the only normal 
       precipitous ideal on $\kappa$. In particular, $L[U]$ satisfies ($*$).
    \item Forcing ($*$), Section \ref{sec:forcnowin}:\\
      $(\alpha)$ We construct a partial order $R(\kappa)$
         forcing that empty wins $\BM_{\leq \omega}(\kappa)$.
         This $R(\kappa)$ does
         not preserve measurability of $\kappa$.\\
      $(\beta)$ We use $R(\kappa)$ to force ($*$) while preserving
         e.g.  supercompactness.
  \end{enumerate}
}

\subsection*{Step 2}

Now we look at the Levy-collapse $Q$ that collapses $\kappa$ to $\theta^+$.

In Section \ref{sec:stillnostrategyafterlevy} we will see:
If in $V[G_Q]$, nonempty wins $\BM_{\leq\omega}(\name S)$
for some $\name S\in \mathcal{E}^\kappa_\theta$, then in $V$
nonempty wins $\BM_{\leq\omega}(\tilde S)$ for some $\tilde S\in \mathcal{E}^\kappa_{\geq \theta}$.

So if we start with $V$ satisfying ($*$) of Step 1, then $Q$ forces:
  \begin{itemize}
    \item Nonempty wins $\PD_{<\theta}(E^\kappa_\theta)$ (by
      \ref{thm:old}.4). Actually nonempty wins $\PD_{<\theta}(S)$
      for all $S\in U$, and $E^\kappa_\theta=(E^\kappa_{\geq \theta})^V\in U$.
    \item Nonempty does not win
      $\BM_{\leq \omega}(\name S)$ for any stationary 
      $\name S\subseteq E^\kappa_\theta$.
      Equivalently: Nonempty does not win $\BM_{\leq \omega}(E^\kappa_\theta)$,
      even if nonempty gets the first move.
  \end{itemize}

\section{$U$ is the only normal, precipitous ideal in $L[U]$}\label{sec:inner}

If $V=L$, then there are no normal, precipitous ideals (recall that a
precipitous ideal implies a measurable in an inner model).  Using Kunen's
results on iterated ultrapowers,
it is easy to relativize this to $L[U]$:
\begin{Thm}\label{cor:inner}
  Assume $V=L[U]$, where $U$ is a normal ultrafilter on $\kappa$.
  Then the dual ideal of $U$ is the only normal, 
  precipitous ideal on $\kappa$. 
  
  In particular, $\text{NS}_\kappa$ is nowhere precipitous, and
  empty wins $\BM_{\leq \omega}(S)$ for any stationary $S\subseteq \kappa$.
\end{Thm}

Remark: Much deeper results by Gitik show that e.g.\\
\centerline{$(\star)$ $\kappa$ is measurable and either $E^\kappa_\lambda$ or 
$\text{NS}_\kappa\restriction\text{Reg}$ is precipitous.}
implies more than a measurable (in an inner model) \cite[Sect. 5]{MR1357746},
so ($\star$) fails not only in $L[U]$ but also in any other universe without
``larger inner-model-cardinals''.
However, it is not clear to us whether the same hold e.g. for\\
\centerline{($\star'$) $\kappa$ is measurable and $\text{NS}\restriction S$ is
precipitous for some $S$.}

Back to the proof of Theorem \ref{cor:inner}.

If $\mathcal{I}$ is a normal, precipitous ideal, then $P_\mathcal{I}$ forces
that the generic filter $G$ is a normal, wellfounded $V$-ultrafilter (cf
\cite[22.13]{MR1940513}). So it is enough to show that in any forcing
extension, $U$ is the only  normal wellfounded $V$-ultrafilter
on $\kappa$.  We will do this in Lemma \ref{lem:gurke3}.

If $U\in L[U]$ and $L[U]$ thinks that $U$ is a normal ultrafilter on $\kappa$,
then we call the pair $(L[U],U)$ a $\kappa$-model.

If $D$ is a normal ultrafilter on $\kappa$, and $U=D\cap L[D]$, then $(L[U],U)$
is a $\kappa$-model.

We will use the following results of Kunen \cite{MR0277346}, cited as
Theorem 19.14 and Lemma 19.16 in \cite{MR1940513}:
\begin{Lem}
  \begin{enumerate}
    \item For every ordinal $\kappa$ there is at most one $\kappa$-model.
    \item Assume $\kappa<\lambda$ are ordinals, $(L[U],U)$ is the
      $\kappa$-model and $(L[W],W)$ the $\lambda$-model.
      Then $(L[W],W)$ is an iterated ultrapower of $(L[U],U)$,
      in particular:
      There is an elementary embedding $i: L[U]\to L[W]$ definable
      in $L[U]$ such that $W=i(U)$. 
    \item Assume that 
      \begin{itemize}
        \item $(L[U],U)$ is the $\kappa$ model, 
        \item $A$ is a set of ordinals of size at least $\kappa^+$,
        \item $\theta$ is a cardinal such that $A\cup \{U\}\subset L_\theta[U]$, and
        \item $X\subseteq \kappa$ is in $L[U]$.
      \end{itemize}
      Then there is a formula $\varphi$, ordinals  $\alpha_i<\kappa$ and
      $\gamma_i\in A$ such that in $L_\theta[U]$,
      $X$ is defined by $\varphi(X,\alpha_1,\dots,\alpha_n,\gamma_1,\dots,\gamma_m,U)$.
  \end{enumerate}
\end{Lem}
(That means that in $L[U]$ there is exactly one $y$ satisfying
$\varphi(y,\alpha_1,\dots)$, and $y=X$.)

\begin{Lem}\label{lem:gurke3}
  Assume $V=L[U]$, where $U$ is a normal ultrafilter on $\kappa$.
  Let $V'$ be a forcing extension of $V$, and $G\in V'$ a normal,
  wellfounded $V$-ultrafilter on $\kappa$. Then $G=U$.
\end{Lem}

\begin{proof}
  In $V'$, let $j:V\to \ult_G(V)$ be elementary. 
  Set $\lambda=j(\kappa)>\kappa$ and $W=j[U]$.
  So $\ult_G(V)$ is the $\lambda$-model $L[W]$.

  In $V$, we can define a function $J:\ON\to\ON$ such that
  in $V'$, $J(\alpha)$ is a cardinal greater than
  $(\alpha^\kappa)^{+V'}$.
  (After all, $V'$ is just a forcing extension of $V$.)
  So $J(\alpha)$ is greater than both $i(\alpha)$ and $j(\alpha)$.
  In $V$, let $\mathcal{C}$ 
  be the class of ordinals that are $\omega$-limits of 
  iterations of $F$, i.e. $\alpha\in \mathcal{C}$ if 
  $\alpha=\sup(\alpha_0,F(\alpha_0),F(F(\alpha_0)),\dots)$. 
  Then $i(\alpha)=j(\alpha)=\alpha$ and $\alpha$ is a cardinal in $V'$.

  In $V'$, pick a set $A$ of $\kappa^+$ many members of $\mathcal{C}$,
  and $\theta \in\mathcal{C}$ such that and $A\cup\{U\}\subseteq L_\theta[U]$. 
  Pick any $X\subseteq \kappa$. Then in $L[U]$, $X$ is defined by
  \[ L_\theta[U]\vDash \varphi(X,\vec\alpha,\vec\gamma,U). \]
  Let $k$ be either $i$ or $j$. Then by elementarity, in $L[W]$ 
  $k(X)$ is the set $Y$ such that
  \[ L_{\theta}[W]\vDash \varphi(Y,\vec\alpha,\vec\gamma,W),\]
  since $W=k(U)$ and $k(\beta)=\beta$ for all $\beta\in \kappa\cup A\cup \{\theta\}$.
  
  Therefore $i(X)=j(X)=Y$. So
  $X\in G$ iff $\kappa\in j(X)=i(X)$ iff $X\in U$, since both $G$ and $U$ are
  normal.
\end{proof}

\section{Forcing empty to win}\label{sec:forcnowin}

As in the last section, we construct a universe with in which empty wins
$\BM_{\leq\omega}(S)$ for every stationary $S\subseteq \kappa$, this time using
forcing.  This shows that  the assumption is also compatible with e.g. $\kappa$
supercompact.

\subsection{The basic forcing}\label{subs:a}

\begin{Asm}\label{asm1}
  $\kappa$ is inaccessible, $2^\kappa=\kappa^+$, and $\lhd$ a wellordering of
  $2^\kappa$ (used for the bookkeeping).
\end{Asm}

We will define the $<\kappa$-support iteration
$(P_\alpha,Q_\alpha)_{\alpha<\kappa^+}$ and show:

\begin{Lem}\label{lem:onestep}
  $P_{\kappa^+}$ forces: Empty has a winning
  strategy for $\BM_{\leq \omega}(\kappa)$ 
  where empty's first move is $\kappa$.
  $P_{\kappa^+}$ is $\kappa^+$-cc and has a dense subforcing
  $P'_{\kappa^+}$ which is $<\kappa$-directed-closed and
  of size $\kappa^+$.
\end{Lem}

We use two basic forcings in the iteration:
\begin{itemize}
  \item If $S\subseteq \kappa$ is stationary, then $\text{Cohen}(S)$
    adds a Cohen subset of $S$.
    Conditions are functions $f:\zeta\rightarrow \{0,1\}$ with
    $\zeta<\kappa$ successor such that
    $\{\xi<\zeta:\, f(\xi)=1\}$ is a subset of $S$.
    $\zeta$ is called height of $f$.
    $\text{Cohen}(S)$ is ordered by inclusion.

    This forcing adds the
    generic set
    $S'=\{\zeta<\kappa:\, (\exists f\in G)f(\zeta)=1\}\subset S$.

  \item If $\lambda\leq\kappa^+$, and $(S_i)_{i<\lambda}$ 
    is a family of stationary sets, then
    $\text{Club}((S_i)_{i<\lambda})$ consists of
    $f:(\zeta\times u)\rightarrow \{0,1\}$, $\zeta<\kappa$ successor,
    $u\subseteq \lambda$, $|u|<\kappa$ 
    such that $\{\xi<\zeta: f(\xi,i)=1\}$ is a closed subset of 
    $S_i$.
    $\zeta$ is called height of $f$, $u$ domain of $f$.
    $\text{Club}((S_i)_{i<\lambda})$ is ordered by inclusion.
\end{itemize}

The following is well known:
\begin{Lem}\label{lem:cohen}
  $\text{Cohen}(S)$ is $<\kappa$-closed and forces that the generic
  Cohen subset $S'\subseteq S$ is stationary.
\end{Lem}

So $\text{Cohen}(S)$ is a well-behaved forcing, adding a generic stationary
subset of $S$. $\text{Club}((S_i)_{i<\lambda})$
adds unbounded closed subsets of each $S_i$.  Other than that
it is not clear why this forcing should e.g.
preserve the regularity of $\kappa$ (and it
will generally not be $\sigma$-closed). 
However, we will shoot clubs only through
complements of Cohen-generics we added previously,
and this will simplify matters considerably.

The $P_\alpha$ will add more and more moves to our winning strategy.

Set $D=\{\delta<\kappa^+:\, \delta\text{ limit}\}$
(for ``destroy''),
$M=(\kappa^+)^{<\omega}$ (for ``moves'').
Find a bijection of $i: M\rightarrow \kappa^+\setminus D$
so that
$s\preceq_M t$ implies $i(s)\leq i(t)$. We identify 
$M$ with its image, i.e. $\kappa^+=D\cup M$.
So for $\alpha\in M$
there is a finite set $\alpha_0<\alpha_1\dots<\alpha_m<\alpha$
of $M$-predecessors (in short: predecessors).
For $\delta\in D$, we can look at all branches through $M\cap \delta$.  Some of
them will be ``new'', i.e. not in any $M\cap \gamma$ for $\gamma<D\cap \delta$.
Let $\lambda_\delta$ be the number of these new branches, i.e.
$0\leq \lambda_\delta\leq 2^\kappa=\kappa^+$.

We define $Q_\alpha$ by induction on $\alpha$, and assume that at stage
$\alpha$ (i.e. after forcing with $P_\alpha$) we have already defined a partial
strategy. Work in $V[G_\alpha]$.

\begin{itemize}
  \item $\alpha\in M$, with the predecessors $0=\alpha_0<\alpha_1\dots<\alpha_m<\alpha$.
    By induction we know that at stage $\alpha_m$
    \begin{itemize}
      \item we dealt with the sequence 
	$x_{\alpha_m}=(\kappa,T_{\alpha_1},S_{\alpha_1},T_{\alpha_2},\dots,
        S_{\alpha_{m-1}},T_{\alpha_{m}})$, which is played to empty's partial
        strategy,
      \item we defined $Q_{\alpha_m}$ to be $\text{Cohen}(T_{\alpha_m})$, adding the generic set $S_{\alpha_m}$,
      \item this $S_{\alpha_m}$ was added to the partial strategy as response to $x_{\alpha_m}$.
    \end{itemize}
    Now (using some simple bookkeeping) we pick a stationary
    $T_{\alpha}\subset S_{\alpha_m}$ such that the partial strategy 
    is not already defined on $x_{\alpha}={x_{\alpha_m}}^\frown ( S_{\alpha_m},T_\alpha)$,
    and set $Q_{\alpha}=\text{Cohen}(T_{\alpha})$,
    and add the $Q_\alpha$-generic
    $S_{\alpha}\in V[G_{\alpha+1}]$ to the partial strategy as response to $x_{\alpha}$.

  \item $\alpha\in D$.
    In $V$, there are $0\leq\lambda_\alpha\leq\kappa^+$ many new branches
    $b_i$. (All old branches have already been dealt with in the previous
    $D$-stages.)
    For each new branch $b_i=(\alpha^i_0<\alpha^i_1<\dots)$, we set 
    $S^i=\bigcap_{n\in\omega} S_{\alpha^i_n}$,
    and we set
    $Q_\alpha=\text{Club}((\kappa\setminus S^i)_{i\in\lambda_\alpha})$.
\end{itemize}

So empty always responds to nonempty's move $T$ with a Cohen subset of $T$, and
the intersection of an $\omega$-sequence of moves according to the strategy is
made non-stationary.

We will show:
\begin{Lem}\label{lem:nonstat}
  $P_{\kappa^+}$ does not add any new countable sequences of ordinals,
  forces that $\kappa$ is regular and that the $Q_\alpha$-generic
  $S_\alpha$ (i.e. empty's move) is stationary for all $\alpha\in M$.
\end{Lem}

We will prove this Lemma later.  Then the rest follows easily:
\begin{Lem}
  $P_{\kappa^+}$ forces that
  empty wins $\BM_{\leq \omega}(\kappa)$, using $\kappa$ as first move.
\end{Lem}

\begin{proof}
  At the final limit stage, $P_{\kappa^+}$ does not add any new
  subsets of $\kappa$,
  nor any countable sequences of such subsets. So there are
  only $\kappa^+$ many names for countable sequences
  $x=(\kappa,T'_1,S'_1,T'_2,S'_2,T'_3,\dots)$.
  Our bookkeeping has to make sure that
  for every initial segment (if it consists of valid moves 
  and uses the partial strategy 
  so far) there has to be a response in the strategy. 

  Then $x\restriction 2n$ corresponds to an element of $M$ 
  for every $n$, and $x$ defines a branch $b$ through $M$.
  $b\in V$, since $P_{\kappa^+}$ does not add new countable sequences
  of ordinals.

  Let $\alpha\in D$ be minimal so that $x\restriction 2n<\alpha$
  for all $n$. 
  Then in the $D$-stage $\alpha$,
  the stationarity of $\bigcap_{n\in\omega} S'_n$ was
  destroyed, i.e. empty wins the run $x$.
\end{proof}

We now define the dense subset of $P_\alpha$:
\begin{Def} $p\in P'_\alpha$ if $p\in P_\alpha$ and
  there are (in $V$) a successor ordinal
  $\epsilon(p)<\kappa$, $(f_\alpha)_{\alpha\in\dom(p)}$ 
  and $(u_\alpha)_{\alpha\in\dom(p)\cap D}$ such that:
  \begin{itemize}
    \item If $\alpha\in M$, then $f_\alpha: \epsilon(p)\rightarrow \{0,1\}$.
    \item If $\alpha\in D$, then $u_\alpha\subseteq \lambda_\alpha$, 
      $|u_\alpha|<\kappa$, and 
      $f_\alpha: \epsilon(p)\times u_\alpha\rightarrow \{0,1\}$.
    \item Moreover, for $\alpha\in D$, $u_\alpha$
      consists exactly of the new branches through $\dom(p)\cap\alpha\cap M$. 
    \item $p\restriction \alpha\Vdash p(\alpha)=f_\alpha$. 
  \end{itemize}
\end{Def}

So a $p\in P'_\alpha$ corresponds to a ``rectangular'' matrix
with entries in $\{0,1\}$.
Of course only some of these matrices are conditions of $P_\alpha$ 
and therefore in $P'_\alpha$.

\begin{Lem}
  \begin{enumerate}
    \item $P'_\alpha$ is ordered by extension. (I.e. if $p,q\in P'_\alpha$,
      then $q\leq p$ iff $q$ (as Matrix) extends $p$.)
    \item $P'_\alpha\subseteq P_\alpha$ is a dense subset.
    \item $P'_\alpha$ is $<\kappa$-directed-closed, in particular
      $P_\alpha$ does not add any new sequences of length $<\kappa$
      nor does it destroy stationarity of any subset of $\kappa$.
  \end{enumerate}
\end{Lem}

\begin{proof}
  (1) should be clear.

  (3) Assume all $p_i$ are pairwise compatible.
  We construct a condition $q$ by putting an additional row on top of
  $\bigcup p_i$ (and filling up at indices where new branches
  might have to be added). So we set
  \begin{itemize}
    \item $\dom(q)=\bigcup\dom(p_i)$.
    \item $\epsilon(q)=\bigcup \epsilon(p_i)+1$.
    \item For $\alpha\in \dom(q)\cap M$, we put $0$ on top, i.e.
      $q_\alpha(\epsilon(q)-1)=0$.
    \item For $\alpha\in \dom(q)\cap D$,
      and $i\in \bigcup \dom(p_i(\alpha))$, set
      $q_\alpha(\epsilon(q)-1,i)=1$.
    \item For $\alpha\in \dom(q)\cap D$, if 
      $i$ is a new branch through $M\cap \dom(q)\cap \alpha$
      and not in $\bigcup \dom(p_i(\alpha))$,
      set $q_\alpha(\xi,i)=0$ for all $\xi<\epsilon(q)$.
  \end{itemize}
  Why can we do that? If $\alpha\in M$,
  whether the bookkeeping says that $\epsilon(q)-1\in T_\alpha$ or not,
  we can of course always choose to not put it into $S_\alpha$
  (i.e. set $q_\alpha(\epsilon(q)-1)=0$).
  Then for $\alpha\in D$, 
  $\epsilon(q)-1$ will definitely not be in the intersection along the 
  branch $i$,
  so we can put it into the complement.
  
  (2) By induction on $\alpha$. Assume $p\in P_\alpha$.

  $\alpha=\beta+1$ is a successor. We know that $P_\beta$ does not add any
  new $<\kappa$ sequences of ordinals, so we can strengthen 
  $p\restriction\beta$ to a $q\in P'_\beta$ which decides $f=p(\beta)\in V$. 
  Without loss of generality $\epsilon(q)\geq \text{height}(f)$,
  and we can enlarge $f$ up to $\epsilon(q)$
  by adding values $0$ (note that $\text{height}(f)<\kappa$ is a successor, so
  we do not get problems with closedness when adding $0$).
  And again, we also add values for the required ``new branches'' if
  necessary.
   
  If $\alpha$ is a limit of cofinality $\geq\kappa$, then $p\in P_\beta$
  for some $\beta<\alpha$, so there is nothing to do.

  Let $\alpha$ be a limit of cofinality $<\kappa$, i.e.
  $(\alpha_i)_{i\in\lambda}$ is an increasing cofinal sequence in $\alpha$, $\lambda<\kappa$.
  Using (2), define a sequence $p_i\in P'_{\alpha_i}$ such that $p_i<p_j\wedge
  p\restriction \alpha_i$ for all $j<i$, then use (3).
\end{proof}

How does the quotient forcing $P^\alpha_{\kappa^+}$ 
(i.e. $P_{\kappa^+}/G_\alpha$) behave compared to $P_{\kappa^+}$?
\begin{itemize}
  \item
    Assume $\alpha\in D$. In $V[G_\alpha]$, $Q_\alpha$
    shoots a club through the complement of the (probably) stationary set
    $\bigcap_{i\in\omega} S^i$. In particular, $Q_\alpha$
    cannot have a $<\kappa$-closed
    subset.
  \item 
    Nevertheless, $P_\alpha\ast Q_\alpha$ has a $<\kappa$-closed
    subset (and preserves stationarity).
  \item 
    So if we factor $P_{\kappa^+}$ at some $\alpha\in D$, the remaining
    $P^\alpha_{\kappa^+}$ will look very different from $P_{\kappa^+}$.
  \item
    However, if we factor $P_{\kappa^+}$ at $\alpha\in M$,
    $P^\alpha_{\kappa^+}$ will be more or less the same as
    $P^\alpha_{\kappa^+}$
    (just with a slightly different bookkeeping).
\end{itemize}
  
In particular, we get:
\begin{Lem}
  If $\alpha\in M$, then the quotient $P^\alpha_{\kappa^+}$ 
  will have a dense $<\kappa$-closed subset 
  (and therefore it will not collapse stationary sets).
\end{Lem}
(The proof is the same as for the last lemma.)

Note that for this result it was necessary to collapse the new branches as soon as
they appear. If we wait with that, then (looking at the rest of the forcing
from some stage $\alpha\in M$) we shoot clubs through stationary sets that
already exist in the ground model, and things get more complicated.

Now we can easily prove lemma \ref{lem:nonstat}:
\begin{proof}[Proof of lemma \ref{lem:nonstat}]
  In stage $\alpha\in M$, nonempty's previous move $S_{\alpha_m}$
  is still stationary (by induction), 
  the bookkeeping chooses a stationary
  subset $T_{\alpha_m}$ of this move, and we
  add $S_\alpha$ as Cohen-generic subset of $T_{\alpha_m}$.
  So according to lemma \ref{lem:cohen}, $S_\alpha$ is stationary
  at stage $\alpha+1$, i.e. in $V[G_{\alpha+1}]$. But since
  $\alpha+1\in M$, the
  rest of the forcing, $P^{\alpha+1}_{\kappa^+}$,
  is $<\kappa$-closed and does not destroy stationarity of $S_\alpha$.
\end{proof}

\subsection{Preserving Measurability}\label{subsec:forcmeas}

We can use the following theorem of Laver \cite{MR0472529},
generalizing an idea of Silver:
If $\kappa$ is supercompact,
then there is a forcing extension in which $\kappa$ is supercompact
and every $<\kappa$-directed closed forcing preserves the supercompactness.
Note that we can also get $2^\kappa=\kappa^+$ which such a forcing.

\begin{Cor}
  If $\kappa$ is supercompact, we can force that 
  $\kappa$ remains supercompact and that empty wins $\BM_{\leq\omega}(S)$
  for all stationary $S\subseteq \kappa$.
\end{Cor}

Remark: It is possible, but not obvious that we can also start with $\kappa$
just measurable and preserve measurability. It is at least likely that 
it is enough to start with strong to get measurable. Much has been published
on such constructions, starting with Silver's proof for violating GCH at
a measurable (as outlined in \cite[21.4]{MR1940513}). 

\section{The Levy collapse}\label{sec:stillnostrategyafterlevy}

We show that after collapsing $\kappa$ to $\theta^+$, nonempty still has no
winning strategy in $\BM$.

Assume that $\kappa$ is inaccessible, $\theta<\kappa$ regular, and 
let $Q = \text{Levy}(\theta, <\kappa)$ be 
the Levy collapse of $\kappa$ to $\theta^+$:
A condition $q \in Q$ is a function defined on a subset of $\kappa \times \theta$, such that
$|\dom(q)| < \theta$ and $q(\alpha, \xi) < \alpha$ for $\alpha>1, (\alpha, \xi) \in \dom(q)$
and $q(\alpha, \xi)=0$ for $\alpha\in\{0,1\}$.

Given $\alpha < \kappa$, define
$Q_\alpha = \{q:\, \dom(q) \subseteq \alpha \times \theta\}$ and
$\pi_\alpha: Q \to Q_\alpha$ by
$q \mapsto q \upharpoonright (\alpha \times \theta)$.

The following is well known (see e.g. \cite[15.22]{MR1940513} for a proof):
\begin{Lem}\label{lem:basiclevy}
  \begin{itemize}
    \item $Q$ is $\kappa$-cc and $<\theta$-closed.
    \item In particular, $Q$ preserves stationarity of subsets of $\kappa$:
      \\
      If $p$ forces that $\name{C}\subseteq \kappa$ is club, then there is a 
      $C'\subseteq \kappa$ club and a $q\leq p$
      forcing that $C'\subseteq \name{C}$.
    \item If $q\Vdash p\in G$, then $q\leq p$ (i.e. $\leq^*$ is the same as $\leq$).
  \end{itemize}
\end{Lem}

We will use the following simple consequence of Fodor's lemma (similar to a $\Delta$-system lemma):

\begin{lemma} \label{lemma100}
  Assume that $p\in Q$ and
  $S \in\mathcal{E}^\kappa_{\geq \theta}$.
  If $\{q_\alpha \mid \alpha \in S\}$ is a sequence of conditions in $Q$, $q_\alpha<p$,
  then there is a $\beta<\kappa$, a $q\in Q_\beta$ and
  a stationary $S' \subseteq S$, such that $q\leq p$ and
  $\pi_\alpha(q_\alpha)=q$ for all $\alpha\in S'$.
\end{lemma}

\begin{proof}
  For $q\in Q$ set $\dom^\kappa(q)=\{\alpha\in\kappa:\, (\exists
  \zeta\in\theta)\, (\alpha,\zeta)\in\dom(q)\}$.
  For $\alpha\in S$ set $f(\alpha) = \sup(\dom^\kappa(q_\alpha) \cap \alpha)$.
  $f$ is regressive,
  since $|\dom^\kappa(q_\alpha)| < \theta$ and $\cf(\alpha) \ge \theta$.
  By the pressing down lemma there is a $\beta<\kappa$ such that
  $T=f^{-1}(\beta)\subseteq S$ is stationary.

  For $\alpha\in T$, set $h(\alpha) = \pi_{\beta+1}(q_\alpha)$. 
  The range of $h$ is of size at most $|\beta\times \theta|^{<\theta} < \kappa$. 
  So  there is a stationary $S' \subseteq T$ such that $h$ is constant on $S'$,
  say $q$.
  If $\alpha\in S'$, then $\sup(\dom^\kappa(q_\alpha)) \cap \alpha=\beta$,
  therefore $\pi_\alpha(q_\alpha)=\pi_{\beta+1}(q_\alpha)=q$.

  Pick $\alpha\in S'$ such that $\alpha>\sup(\dom^\kappa(p))$. $q_\alpha\leq p$, so
  $q=\pi_{\alpha}(q_\alpha)\leq \pi_{\alpha}(p)=p$.
\end{proof}

\begin{lemma} Assume that
\begin{itemize}
  \item	$\kappa$ is strongly inaccessible, $\theta<\kappa$ regular, $\mu\leq \theta$.
  \item	$Q = \text{Levy}(\theta, <\kappa)$,
  \item $\name S$ is a $Q$-name for an element of $\mathcal{E}^\kappa_{\theta}$,
  \item	$\tilde p\in Q$ forces that $\name{F}$
	is a winning strategy of nonempty in $\BM_{<\mu}(\name S)$.
\end{itemize}
Then in $V$, nonempty wins $\BM_{<\mu}(\tilde S)$ for some 
$\tilde S\in E^\kappa_{\geq\theta}$. 

If $\name S$ is a standard name for $T\in (E^\kappa_{\geq\theta})^V$, 
then we can set $S=T$.
\end{lemma}

\begin{proof}
  First assume that $\name S$ is a standard name.

  For a run of $\BM_{<\mu}(S)$, we let 
  $A_\varepsilon$ and $B_\varepsilon$ denote the $\varepsilon$th moves of empty and nonempty.
  We will construct by induction 
  on $\varepsilon < \mu$ 
  a strategy for empty, including not only the moves $B_\varepsilon$, but also
  $Q$-names $\name{A}_\varepsilon', \name{B}_\varepsilon'$, and $Q$-conditions
  $p_\varepsilon, \langle p_\alpha^\varepsilon \mid \alpha \in B_\varepsilon \rangle$, such that
  the following holds:
  \begin{itemize}
  \item	$p_\varepsilon \le p_\xi$ and  $p^\varepsilon_\alpha \le p^\xi_\alpha$ for $\xi < \varepsilon$.
  \item	$p_\varepsilon$ forces that $(\name{A}_\xi', \name{B}_\xi')_{\xi \le \varepsilon}$
  	is an initial segment of a run of $\BM_{<\mu}(\name S)$ in which
  	nonempty uses the strategy $\name{F}$.
  \item	$p_\varepsilon \Vdash \name{A}_\varepsilon' \subseteq A_\varepsilon$.
  \item	For $\alpha \in B_\varepsilon$, 
  	$\pi_\alpha(p^\varepsilon_\alpha) = p_\varepsilon$ (in particular $p^\varepsilon_\alpha\le p_\varepsilon$), and
  	$p^\varepsilon_\alpha \Vdash$ ``$\alpha \in \name{B}_\varepsilon'$''.
  
  \end{itemize}

  Assume that we have already constructed these objects for all $\xi < \varepsilon$.
  
  In limit stages $\varepsilon$, we first have to make sure that 
  $\bigcap_{\xi < \varepsilon} B_\xi$ is stationary (otherwise nonempty has already lost).
  Pick a $q$ stronger than each $p_\xi$ for $\xi < \varepsilon$. (This is possible since $Q$ is $<\theta$-closed.)
  Then $q$ forces that $\bigcap_{\xi < \varepsilon} B_\xi=\bigcap_{\xi < \varepsilon} A_\xi\supseteq \bigcap_{\xi < \varepsilon}  \name{A}'_\xi$
  and that $(\name{A}_\xi', \name{B}_\xi')_{\xi \le \varepsilon}$
  is a valid initial segment of a run where nonempty uses the strategy, in particular
  $\bigcap_{\xi < \varepsilon}  \name{A}'_\xi$ is stationary.
  
  So now $\varepsilon$ can be a successor or a limit, and empty 
  plays the stationary set $A_\varepsilon \subseteq \bigcap_{\xi < \varepsilon} B_\xi$.
  (That implies that $p_\alpha^\xi$ is defined for all $\alpha\in A_\varepsilon$ and
  $\xi < \varepsilon$.)
  \begin{itemize}
    \item
      Define the $\varepsilon$th move of empty in $V[G_Q]$ to be 
      \begin{displaymath}
         \name{A}_\varepsilon' = \{\alpha \in A_\varepsilon:\,  (\forall \xi < \varepsilon)\,
  	p_\alpha^\xi \in G_Q \},
      \end{displaymath}
      and pick $\tilde p_\varepsilon \le p_\xi$ for $\xi < \varepsilon$
      (for $\varepsilon=0$, pick $\tilde p_0=\tilde p$).
  
      $\tilde p_\varepsilon$ forces that
      $\name{A}_\varepsilon'\subseteq \bigcap_{\xi < \varepsilon}  \name{B}'_\xi$, since 
      $p_\alpha^\xi$ forces that $\alpha\in \name{B}'_\xi$.
      $\tilde p_\varepsilon$ also forces that $\name{A}_\varepsilon'$ is stationary:
      \\
      Otherwise there is a $C\subseteq \kappa$ club and a $q\leq \tilde p_\varepsilon$
      forcing that 
      $C \cap \name{A}_\varepsilon$ is empty (cf \ref{lem:basiclevy}). $q \in Q_\beta$
      for some $\beta < \kappa$. Pick $\alpha\in (C\cap A_\varepsilon)\setminus (\beta+1)$.
      For $\xi < \varepsilon$, $\pi_\alpha(p^\xi_\alpha)=p_\xi\geq q$,
      and $q\in Q_\beta$,
      so $q$ and $p^\xi_\alpha$ are compatible. Moreover, the conditions
      $(q\cup p^\xi_\alpha)_{\xi\in\varepsilon}$ are decreasing, so there
      is a common lower bound $q'$ forcing that $p^\xi_\alpha\in G_Q$ for all $\xi$,
      i.e. that $\alpha\in \name{A}_\varepsilon'$, a contradiction.
    \item 
      Given $\name{A}_\varepsilon'$, we define $\name{B}_\varepsilon'$ as the response
      according to the strategy $\name{F}$.
    \item
      Now we show how to obtain the next move of nonempty, $B_\varepsilon$,
      (in the ground model), as well as $p^\varepsilon_\alpha$ for $\alpha\in B_\varepsilon$.
      $B_\varepsilon$ of course has to be a subset of the stationary set $S$
      defined by
      \begin{displaymath}
        S = \{\alpha \in A_\varepsilon \mid \tilde p_\varepsilon \not\Vdash \alpha \notin 
        \name{B}_\varepsilon' \}.
      \end{displaymath}
  
      For each $\alpha \in S$, pick some $p^\varepsilon_\alpha \le \tilde p_\varepsilon$ forcing
      that $\alpha \in \name{B}_\varepsilon'$.
      By the definition of $\name{A}_\varepsilon'$ and since 
      $\tilde p_\varepsilon\Vdash \name{B}_\varepsilon' \subseteq \name{A}_\varepsilon'$ , we get
      \begin{displaymath} 
        p^\varepsilon_\alpha \Vdash (\forall \xi < \varepsilon)\, p_\alpha^\xi \in G_Q,
      \end{displaymath}
      which means that for $\alpha \in S$ and $\xi < \varepsilon$, 
      $p^\varepsilon_\alpha \le p^\xi_\alpha$.
  
      Now we apply lemma \ref{lemma100} (for $p=\tilde p_\varepsilon$). This gives us 
      $S'\subseteq S$ and $q\leq \tilde p_\varepsilon$. We set
      $B_\varepsilon=S'$ and $p_\varepsilon=q$.
  \end{itemize}

  If $\name S$ is not a standard name, set 
  \[
    S^0=\{ \alpha\in E^\kappa_{\geq \theta}:\, \tilde p\not\Vdash
    \alpha\notin\name S  \}
  \]
  As above, for each $\alpha\in S_0$,
  pick a $\tilde p^{-1}_\alpha\leq  \tilde p$ forcing that $\alpha\in\name S$,
  and choose a stationary $\tilde S\subseteq S^0$ according to Lemma
  \ref{lemma100}.
  Now repeat the proof, starting the sequence $(p_\varepsilon)$ and
  $(p^\varepsilon_\alpha)$ already at $\varepsilon=-1$.
\end{proof}

\bibliographystyle{amsplain}
\bibliography{896}

\providecommand{\bysame}{\leavevmode\hbox to3em{\hrulefill}\thinspace}
\providecommand{\MR}{\relax\ifhmode\unskip\space\fi MR }
\providecommand{\MRhref}[2]{%
  \href{http://www.ams.org/mathscinet-getitem?mr=#1}{#2}
}
\providecommand{\href}[2]{#2}
\begin{thebibliography}{10}

\bibitem{MR716633}
Matthew Foreman, \emph{Games played on {B}oolean algebras}, J. Symbolic Logic
  \textbf{48} (1983), no.~3, 714--723. \MR{MR716633 (85h:03064)}

\bibitem{MR0485391}
F.~Galvin, T.~Jech, and M.~Magidor, \emph{An ideal game}, J. Symbolic Logic
  \textbf{43} (1978), no.~2, 284--292. \MR{MR0485391 (58 \#5237)}

\bibitem{MR1357746}
Moti Gitik, \emph{Some results on the nonstationary ideal}, Israel J. Math.
  \textbf{92} (1995), no.~1-3, 61--112. \MR{MR1357746 (96k:03108)}

\bibitem{MR1780068}
Moti Gitik and Saharon Shelah, \emph{Cardinal preserving ideals}, J. Symbolic
  Logic \textbf{64} (1999), no.~4, 1527--1551. \MR{MR1780068 (2002a:03100)}

\bibitem{MR1971240}
Tapani Hyttinen, Saharon Shelah, and Jouko V{\"a}{\"a}n{\"a}nen, \emph{More on
  the {E}hrenfeucht-{F}ra\"\i ss\'e game of length {$\omega\sb 1$}}, Fund.
  Math. \textbf{175} (2002), no.~1, 79--96. \MR{MR1971240 (2004b:03046)}

\bibitem{MR560220}
T.~Jech, M.~Magidor, W.~Mitchell, and K.~Prikry, \emph{Precipitous ideals}, J.
  Symbolic Logic \textbf{45} (1980), no.~1, 1--8. \MR{MR560220 (81h:03097)}

\bibitem{MR519808}
Thomas Jech, \emph{A game theoretic property of {B}oolean algebras}, Logic
  Colloquium '77 (Proc. Conf., Wroc\l aw, 1977), Stud. Logic Foundations Math.,
  vol.~96, North-Holland, Amsterdam, 1978, pp.~135--144. \MR{MR519808
  (80c:90184)}

\bibitem{MR1940513}
\bysame, \emph{Set theory}, Springer Monographs in Mathematics,
  Springer-Verlag, Berlin, 2003, The third millennium edition, revised and
  expanded. \MR{MR1940513 (2004g:03071)}

\bibitem{MR0505504}
Thomas Jech and Karel Prikry, \emph{On ideals of sets and the power set
  operation}, Bull. Amer. Math. Soc. \textbf{82} (1976), no.~4, 593--595.
  \MR{MR0505504 (58 \#21618)}

\bibitem{MR739910}
Thomas~J. Jech, \emph{More game-theoretic properties of {B}oolean algebras},
  Ann. Pure Appl. Logic \textbf{26} (1984), no.~1, 11--29. \MR{MR739910
  (85j:03110)}

\bibitem{MR0739911}
\bysame, \emph{Some properties of {$\kappa $}-complete ideals defined in terms
  of infinite games}, Ann. Pure Appl. Logic \textbf{26} (1984), no.~1, 31--45.
  \MR{MR739911 (85h:03057)}

\bibitem{MR0277346}
Kenneth Kunen, \emph{Some applications of iterated ultrapowers in set theory},
  Ann. Math. Logic \textbf{1} (1970), 179--227. \MR{MR0277346 (43 \#3080)}

\bibitem{MR0472529}
Richard Laver, \emph{Making the supercompactness of {$\kappa $} indestructible
  under {$\kappa $}-directed closed forcing}, Israel J. Math. \textbf{29}
  (1978), no.~4, 385--388. \MR{MR0472529 (57 \#12226)}

\bibitem{MR1469093}
Kecheng Liu and Saharon Shelah, \emph{Cofinalities of elementary substructures
  of structures on {$\aleph\sb \omega$}}, Israel J. Math. \textbf{99} (1997),
  189--205. \MR{MR1469093 (98m:03100)}

\bibitem{MR1191613}
Alan Mekler, Saharon Shelah, and Jouko V{\"a}{\"a}n{\"a}nen, \emph{The
  {E}hrenfeucht-{F}ra\"\i ss\'e-game of length {$\omega\sb 1$}}, Trans. Amer.
  Math. Soc. \textbf{339} (1993), no.~2, 567--580. \MR{MR1191613 (94a:03058)}

\bibitem{MR1420262}
Saharon Shelah, \emph{Large normal ideals concentrating on a fixed small
  cardinality}, Arch. Math. Logic \textbf{35} (1996), no.~5-6, 341--347.
  \MR{MR1420262 (97m:03078)}

\bibitem{MR846604}
Boban Veli{\v{c}}kovi{\'c}, \emph{Playful {B}oolean algebras}, Trans. Amer.
  Math. Soc. \textbf{296} (1986), no.~2, 727--740. \MR{MR846604 (88a:06017)}

\end{thebibliography}

\end{document}